\documentclass[11pt,a4paper]{article} %{amsart} 
\usepackage[final]{optional}
\usepackage[UKenglish]{babel}
\usepackage{bbm}
\usepackage{mathrsfs}

\usepackage{amsfonts, amsmath, wasysym}
\usepackage{amssymb,amsthm,
paralist
}
\usepackage{
latexsym,
}
\usepackage[usenames]{color}

\usepackage{url}

\definecolor{darkgreen}{rgb}{0,0.5,0}
\definecolor{darkred}{rgb}{0.7,0,0}
\usepackage[colorlinks, 
citecolor=darkgreen, linkcolor=darkred
]{hyperref}
\usepackage{esint}
\usepackage{bibgerm}
\usepackage[normalem]{ulem}

\theoremstyle{plain}
\newtheorem{lemma}{Lemma}[section]
\newtheorem{thm}[lemma]{Theorem}

\theoremstyle{definition}

\newtheorem{rmk}[lemma]{Remark}
\def\blbox{\quad \vrule height7.5pt width4.17pt depth0pt}
\newcommand{\cmt}[1]{\opt{draft}{\textcolor[rgb]{0.5,0,0}{
$\LHD$ #1 $\RHD$\marginpar{\blbox}}}}

\numberwithin{equation}{section}

\newcommand{\pl}[2]{{\frac{\partial #1}{\partial #2}}}
\newcommand{\ka}{\kappa}

\newcommand{\vph}{\varphi}

\newcommand{\R}{\ensuremath{{\mathbb R}}}
\newcommand{\N}{\ensuremath{{\mathbb N}}}

\newcommand{\Z}{\ensuremath{{\mathbb Z}}}

\newcommand{\weakto}{\rightharpoonup}

\newcommand{\upto}{\uparrow}
\newcommand{\embed}{\hookrightarrow}

\newcommand{\lap}{\Delta}

\def\blbox{\quad \vrule height7.5pt width4.17pt depth0pt}

\newcommand{\beq}{\begin{equation}}
\newcommand{\eeq}{\end{equation}}
\newcommand{\beqs}{\begin{equation}}
\newcommand{\eeqs}{\end{equation}}
\newcommand{\beqa}{\begin{equation}\begin{aligned}}
\newcommand{\eeqa}{\end{aligned}\end{equation}}
\newcommand{\beqas}{\begin{equation}\begin{aligned}}
\newcommand{\eeqas}{\end{aligned}\end{equation}}
\newcommand{\brmk}{\begin{rmk}}
\newcommand{\ermk}{\end{rmk}}
\newcommand{\partref}[1]{\hbox{(\csname @roman\endcsname{\ref{#1}})}}
\newcommand{\half}{\frac{1}{2}}

\newcommand{\avint}{{\int\!\!\!\!\!\!-}}

\title{{ \bf
{\large A RIGIDITY ESTIMATE FOR MAPS FROM $S^2$ TO $S^2$}\\
{\large VIA THE HARMONIC MAP FLOW}
}
}
\author{Peter M. Topping}
\date{22 September 2020}
\begin{document}
\maketitle

\begin{abstract}
We show how a rigidity estimate introduced in recent work of Bernand-Mantel, Muratov and Simon \cite{muratov} can be derived from the harmonic map flow theory
in \cite{rigidity}.
\end{abstract}

\section{Introduction}

For maps $u\in W^{1,2}(S^2,S^2)$,  consider the harmonic map energy
$$E(u)=\frac12\int|Du|^2.$$
Following Schoen-Uhlenbeck \cite{SU} we know that such a map can be approximated in $W^{1,2}(S^2,S^2)$ by smooth maps $u_i$, whose degree will stabilise at some integer that gives a well-defined notion of degree of $u$. See also \cite[Theorem 6.2]{struwe_book}.
For smooth maps $w$ from a surface to any Riemannian manifold, a simple computation confirms that $E(w)$ is always at least the area of the image of $w$, with equality if $w$ is weakly conformal. This implies that degree $k\in \N$ maps from $S^2$ to itself must have energy at least $4\pi k$, with equality implying that the map is a rational map. A special case of this is that if $u$ is of degree one then $E(u)\geq 4\pi$ with equality if and only if $u$ is a M\"obius map.

%\cmt{add ref to John Wood, or Lemaire? Uses result in book of Springer?}

The following theorem was recently proved by Bernand-Mantel, Muratov and Simon \cite{muratov}. 
\begin{thm}
\label{rig_thm}
There exists a universal constant $C<\infty$ such that for each 
$u\in W^{1,2}(S^2,S^2)$ of degree one, there exists a M\"obius map $v:S^2\to S^2$ such that
$$\int |D(u-v)|^2
%+\int |u-v|^2
\leq C[E(u)-4\pi].$$
\end{thm}

The purpose of this note is to show that the estimate above can be derived from the theory of the harmonic map flow developed in \cite{rigidity}.

\section{The harmonic map flow from surfaces}

The harmonic map flow \cite{ES} is the $L^2$-gradient flow for the harmonic map energy. If we view a flow $u:S^2\times [0,T)\to S^2\embed \R^3$ as taking values in $\R^3$, then the equation in this case can be written 
$$\pl{u}{t}=\tau(u):=(\lap u)^T$$
where $(\lap u)^T$ is the projection of $\lap u$ onto the tangent space of the target, i.e. $(\lap u)^T=\lap u +u|Du|^2$. The energy $E(t):=E(u(t))$ decays according to
\beq
\label{en_decay}
\frac{dE}{dt}=-\|\tau(u)\|_{L^2(S^2)}^2.
\eeq

In 1985 Struwe \cite{struwe_CMH} initiated a theory for the harmonic map flow in the case that the domain is a surface, as it is here. He showed how one can start a flow with smooth or even $W^{1,2}$ initial data, giving a 
global weak solution that is smooth away from finitely many points in space-time, at which bubbling occurs.
At each finite-time singularity, concentrated energy at the singular points is thrown away by taking a weak limit, and the flow restarted. Thus, the energy drops down at each singular time by at least the minimum energy of one bubble, i.e. of a nonconstant harmonic map from $S^2$ to the target. 
However, the map $t\mapsto u(t)$ is continuous into $L^2$, even across singular times. 
(Later \cite{revbub}, different continuations through certain singularities were constructed that did not require a drop in energy, but we will not need them here.)

\cmt{ref to struwe book for infinite time smooth convergence - theorem 6.6 in some editions}

Concerning the asymptotics at infinite time, there exists a sequence $t_i\to\infty$ such that  maps $u(t_i)$ converge smoothly to a limiting harmonic map $u_\infty$ away from finitely many bubble points.
(See, for example, \cite{struwe_book}.)
In general, even in the absence of bubbling,  convergence of the form $u(t)\to u_\infty$ as $t\to\infty$ may fail, even in $L^1$, see \cite{Topping1996, rigidity}. However, a theory was developed in \cite{rigidity, repulsion} concerning such uniform convergence for maps from $S^2$ to itself. A key lemma from that work, which will be useful to us now, gives a relationship between 
the tension field of a map $u$ and its excess energy. It differed from previous estimates of `Lojasiewicz-Simon' type \cite{Simon1983} in that it 
could handle singular objects. In particular the map $u$ below is not asked to be $W^{1,2}$ close to a harmonic map.

\begin{lemma}[{\cite[Lemma 1]{rigidity}}]
\label{key_lemma}
There exist universal constants $\epsilon_0>0$ and $\kappa>0$ such that if a smooth degree $k\in\Z$ map $u:S^2\to S^2$ satisfies 
$E(u)-4\pi |k|<\epsilon_0$,
then 
$$E(u)-4\pi |k|\leq \kappa^2 \|\tau(u)\|_{L^2}^2.$$
\end{lemma}
As is well known, such an estimate gives control on the gradient flow.
Indeed, given a smooth solution $u:S^2\times [0,T]\to S^2$ that is of degree $k\in \Z$, and satisfies $E(u)-4\pi|k|< \epsilon_0$ at time $t=0$ (and therefore also for later times)
we can compute using \eqref{en_decay} and then Lemma \ref{key_lemma} that
$$
-\frac{d}{dt}\left[E(u)-4\pi |k|\right]^\half
=\half \left[E(u)-4\pi |k|\right]^{-\half}\|\tau(u)\|_{L^2}^2
\geq \frac{1}{2\ka}\|\tau(u)\|_{L^2}
$$
and integrating from $t=s\in [0,T)$ to $t=T$ gives
\beqa
\label{tau_int}
\int_s^T \|\tau(u)\|_{L^2}dt
&\leq 2\ka \left([E(u(s))-4\pi |k|]^{\frac12}-[E(u(T))-4\pi |k|]^{\frac12}\right)\\
&\leq 2\ka\left[E(u(s))-4\pi|k|\right]^{\frac12}.
\eeqa
Since $\frac{\partial u}{\partial t}=\tau(u)$, the flow then cannot move far in $L^2$:
\beq
\label{L2_control}
\|u(T)-u(s)\|_{L^2}\leq 2\ka\left[E(u(s))-4\pi|k|\right]^{\frac12}.
\eeq

\section{{Proof of Theorem \ref{rig_thm}}}

\cmt{Here's a bit more detail for the approximation claim below.
Suppose we can prove the estimate for smooth $u$, with a given $C=C_0$.
For a given $u\in W^{1,2}$, take a smooth approximating sequence $u_i$.
Then there exist M\"obius maps $v_i$ such that  
$$\|u_i-v_i\|_{W^{1,2}}^2\leq C_0[E(u_i)-4\pi].$$
Therefore 
$$\|u-v_i\|_{W^{1,2}}^2
\leq 2\|u-u_i\|_{W^{1,2}}^2+2\|u_i-v_i\|_{W^{1,2}}^2
\leq 2\|u-u_i\|_{W^{1,2}}^2+C_0[E(u_i)-4\pi].$$
We can then fix an $i$ sufficiently large so that
$\|u-u_i\|_{W^{1,2}}^2\leq E(u)-4\pi$ and so that 
$E(u_i)-4\pi\leq 2[E(u)-4\pi]$. For this fixed $i$ we then have
$$\|u-v_i\|_{W^{1,2}}^2
\leq 2(1+C_0)[E(u)-4\pi].$$
}

\begin{proof}
First, if $E(u)=4\pi$, then $u$ must be a M\"obius map and we can choose $v=u$, so from now on we may assume that $E(u)>4\pi$.
By approximation, using the definition of degree, it suffices to prove the result for $u\in C^\infty(S^2,S^2)$.
Since the estimate is invariant under pre-composition by M\"obius maps, it suffices to prove the theorem for $u$ equal to a map $u_0$ with the property that 
$\int_{S^2} u_0 =0\in\mathbb{R}^3$. That this can be achieved follows from  a  topological argument: 
For $a\in B^3$, let $\varphi_a:S^2\to S^2$ be the M\"obius map that fixes $\pm\frac{a}{|a|}$, whose differential does not rotate the tangent spaces at $\pm\frac{a}{|a|}$, 
and which when extended to a conformal map $B^3\mapsto B^3$ will send the origin to $a$.
Thus as $a$ approaches some $a_0\in S^2=\partial B^3$, the composition $u\circ \varphi_a$
converges to $u(a_0)$ away from $-a_0$, so $\frac{1}{4\pi}\int u\circ \varphi_a\to u(a_0)$.
Thus the map
$$a\mapsto \frac{1}{4\pi}\int u\circ \varphi_a$$ %\to u(a_0)$$
extends to a continuous map $\Phi$ from $\overline {B^3}$ to itself that agrees with the degree one map $u$ on the boundary $S^2$.
A topological argument then tells us that $\Phi$ is surjective, since otherwise $\Phi$ could be homotoped to a continuous map $\tilde \Phi$ from $\overline {B^3}$ to $S^2$ that restricts to  $u:S^2\to S^2$. But then $\tilde \Phi$ would provide a homotopy from the degree one map $u$ to a constant map, a contradiction.
%we could construct a retraction from $B^3$ to $S^2$.
In particular, there exists $a\in B^3$ such that $\Phi(a)=0\in B^3$.
We can then set $u_0:=u\circ \varphi_a$ to achieve our objective. 
For a related argument in which one \emph{post}-composes with M\"obius maps
see Li-Yau \cite{LY}.

Let now $\epsilon_0$ be as in the key lemma \ref{key_lemma}. For later use, if necessary we reduce $\epsilon_0>0$ so that 
\beq
\label{ep_ineq}
(1+4\ka^2)\epsilon_0\leq \pi.
\eeq
We may assume that our map $u_0$ satisfies $E(u_0)-4\pi<\epsilon_0$ since otherwise the theorem is vacuously true.

\cmt{Elaboration: the LHS of the estimate we're trying to prove is controlled by $C+C[E(u)-4\pi]$}

To prove the theorem, run the harmonic map flow starting with $u_0$.  
By \eqref{L2_control}, the flow is constrained in how far it can move in $L^2$.
Because of the balancing $\int_{S^2} u_0 =0$, this implies that $\int_{S^2} u(t)$ remains close to the origin, which precludes bubbling both at finite and infinite time.
More precisely, because
$$\frac{d}{dt}\int u=\int \tau(u)$$
we have
$$\bigg|\frac{d}{dt}\int u\bigg|\leq 2\sqrt{\pi}%(4\pi)^\half
\|\tau(u)\|_{L^2},$$
and therefore, integrating from $0$ to $t$ using \eqref{tau_int} we have
\beq
\label{com_ineq}
\bigg|\avint u(t) \bigg|\leq \frac{\ka}{\sqrt{\pi}} \left[E(u_0)-4\pi\right]^{\frac12}.
\eeq
Suppose we develop a singularity at a finite or infinite time
$T\in (0,\infty]$. If we pick $t_i\upto T$ such that 
$u(t_i)\weakto w$
weakly in $W^{1,2}$, then we must have
$E(w)\leq E(u_0)-4\pi$, since the singularity loses at least $4\pi$ of energy (that being the least possible energy of a nonconstant harmonic map from $S^2$ to itself) but also because 
$$\avint u(t_i)\to \overline w:=\avint w$$ 
we have
$$|\overline w|\leq \frac{\ka}{\sqrt{\pi}} \left[E(u_0)-4\pi\right]^{\frac12}.$$
%\leq 4\sqrt{2}\ka \epsilon_0.$$
But the Poincar\'e inequality tells us that
$$\int_{S^2}|w-\overline w|^2\leq E(w),$$
and so integrating the inequality
$$1=|w|^2\leq 2|w-\overline w|^2+2|\overline w|^2,$$
we obtain
\beqa
4\pi &\leq 2\int  |w-\overline w|^2 + 8\pi |\overline w|^2\\
&\leq 2 E(w)+8 \ka^2 \left[E(u_0)-4\pi\right]\\
&\leq (2+8 \ka^2)\left[E(u_0)-4\pi\right]\\
&\leq 2\pi,
\eeqa
by \eqref{ep_ineq}, giving a contradiction.

We deduce that the flow exists for all time and converges smoothly to a M\"obius map $v$. Here we only need convergence at some sequence of times $t_i\to\infty$, although we have convergence as $t\to\infty$ by \cite{rigidity}.

We can also pass the inequality \eqref{com_ineq} to the limit $t\to\infty$ to give
$$|\overline v|\leq \frac{\ka}{\sqrt{\pi}} \left[E(u_0)-4\pi\right]^{\frac12}
\leq \ka\left(\frac{\epsilon_0}{\pi}\right)^\half\leq {\textstyle \half},$$
by \eqref{ep_ineq}.
This estimate prevents $v$ from being too concentrated, and we can deduce that 
\beq
\label{Dv_est}
|Dv|\leq c_0,
\eeq
for some universal $c_0$. This follows by noticing that every M\"obius map can be written as one of the maps $\vph_a$ followed by a rotation of $S^2$. Another way of making this precise is by arguing by contradiction: If such an estimate \eqref{Dv_est} were not true, then we would take a sequence of M\"obius maps $v_i$ with $|\overline{v_i}|\leq \half$
but with $\sup|Dv_i|\to\infty$. After a bubbling analysis (passing to a subsequence) the maps $v_i$ would converge weakly in $W^{1,2}$ to a constant map $v_\infty:S^2\to S^2$ with $|\overline{v_\infty}|\leq \half$, a contradiction.

Returning to \eqref{L2_control}, we find that
$\|u_0-v\|_{L^2}^2\leq 4\ka^2[E(u_0)-4\pi]$. But we can also compute
$$\int|D(u_0-v)|^2 = \int |Du_0|^2 +\int |Dv|^2 -2\int\langle Du_0,Dv\rangle$$
and because $-\Delta v=v|Dv|^2$, we can handle the final term using
\begin{equation}
\begin{aligned}
-2\int\langle Du_0,Dv\rangle &= -2\int u_0 (-\Delta v) = -2\int u_0 v|Dv|^2\\
&= \int |u_0-v|^2|Dv|^2 - 2\int |Dv|^2,
\end{aligned}
\end{equation}
where we have used that $|u_0|=|v|=1$. Combining, we obtain
\begin{equation}
\begin{aligned}
\int|D(u_0-v)|^2 &= \int |Du_0|^2 - \int |Dv|^2 + \int |u_0-v|^2|Dv|^2\\
&\leq 2[E(u_0)-4\pi] + c_0^2\|u_0-v\|_{L^2}^2\\
&\leq C[E(u_0)-4\pi],
\end{aligned}
\end{equation}
for universal $C$.
\end{proof}

\brmk
At the start of the argument we balanced our map to have `centre of mass' at the origin. Without this step we could expect the harmonic map flow to generate a finite-time singularity. Indeed in \cite[Theorem 5.5]{Topping1996} we showed that for arbitrarily small $\epsilon_0>0$, there exists a smooth degree one map $u_0:S^2\to S^2$ with $E(u_0)<4\pi+\epsilon_0$ such that the subsequent harmonic map flow must develop a singularity in finite time.
\ermk

{\sc Mathematics Institute, University of Warwick, Coventry,
CV4 7AL, UK}

\end{document}